\documentclass[a4paper,12pt]{amsart}
\usepackage{amssymb}
\usepackage{ifthen}
\usepackage{graphicx}

\newtheorem{thm}{Theorem}
\newtheorem{cor}{Corollary}
\newtheorem{lem}{Lemma}

\newtheorem{que}{Question}
\newtheorem{rem}{Remark}
\newtheorem{prob}{Problem}
\newtheorem{conj}{Conjecture}
\theoremstyle{definition}
\newtheorem{example}[equation]{Example}

\newcounter {own}
\def\theown {\thesection       .\arabic{own}}

\newenvironment{pf}[1][]{%
 \vskip 3mm
 \noindent
 \ifthenelse{\equal{#1}{}}%
  {{\slshape Proof. }}%
  {{\slshape #1.} }%
 }%
{\qed\bigskip}

\newcounter{alphabet}
\newcounter{tmp}
\newenvironment{Thm}[1][]{\refstepcounter{alphabet}%
\bigskip%
\noindent%
{\bf Theorem \Alph{alphabet}}%
\ifthenelse{\equal{#1}{}}{}{ (#1)}%
{\bf .} \itshape}{\vskip 8pt}

\makeatletter
\newcommand{\Ref}[1]{\@ifundefined{r@#1}{}{\setcounter{tmp}{\ref{#1}}\Alph{tmp}}}
\makeatother

\newcounter{minutes}
\divide\time by 60
\newcounter{hours}
\multiply\time by 60
\addtocounter{minutes}{-\time}

\newenvironment{Lem}[1][]{\refstepcounter{alphabet}%
\bigskip%
\noindent%
{\bf Lemma \Alph{alphabet}}%
{\bf .} \itshape}{\vskip 8pt}

\newcommand{\IC}{{\mathbb C}}
\newcommand{\ID}{{\mathbb D}}

\def\be{\begin{equation}}
\def\ee{\end{equation}}

\newcommand{\bee}{\begin{enumerate}}
\newcommand{\eee}{\end{enumerate}}

\newcommand{\blem}{\begin{lem}}
\newcommand{\elem}{\end{lem}}
\newcommand{\bthm}{\begin{thm}}
\newcommand{\ethm}{\end{thm}}
\newcommand{\bcor}{\begin{cor}}
\newcommand{\ecor}{\end{cor}}
\newcommand{\beg}{\begin{example}}
\newcommand{\eeg}{\end{example}}
\newcommand{\bques}{\begin{que}}
\newcommand{\eques}{\end{que}}
\newcommand{\begs}{\begin{examples}}
\newcommand{\eegs}{\end{examples}}
\newcommand{\bdefe}{\begin{defin}}
\newcommand{\edefe}{\end{defin}}
\newcommand{\bprob}{\begin{prob}}
\newcommand{\eprob}{\end{prob}}
\newcommand{\bei}{\begin{itemize}}
\newcommand{\eei}{\end{itemize}}

\newcommand{\bcon}{\begin{conj}}
\newcommand{\econ}{\end{conj}}
\newcommand{\bcons}{\begin{conjs}}
\newcommand{\econs}{\end{conjs}}
\newcommand{\bprop}{\begin{propo}}
\newcommand{\eprop}{\end{propo}}
\newcommand{\br}{\begin{rem}}
\newcommand{\er}{\end{rem}}
\newcommand{\brs}{\begin{rems}}
\newcommand{\ers}{\end{rems}}
\newcommand{\bo}{\begin{obser}}
\newcommand{\eo}{\end{obser}}
\newcommand{\bos}{\begin{obsers}}
\newcommand{\eos}{\end{obsers}}
\newcommand{\bpf}{\begin{pf}}
\newcommand{\epf}{\end{pf}}
\newcommand{\ba}{\begin{array}}
\newcommand{\ea}{\end{array}}
\newcommand{\beq}{\begin{eqnarray}}
\newcommand{\beqq}{\begin{eqnarray*}}
\newcommand{\eeq}{\end{eqnarray}}
\newcommand{\eeqq}{\end{eqnarray*}}

\newcommand{\ra}{\rightarrow}

\newcommand{\ds}{\displaystyle}

\def\cc{\setcounter{equation}{0}   
\setcounter{figure}{0}\setcounter{table}{0}}

\begin{document}

\bibliographystyle{amsplain}

%

\title[Coefficient problems  on the class $U(\lambda)$]{Coefficient problems on the class $U(\lambda)$}

\thanks{
File:~\jobname .tex,
          printed: \number\day-\number\month-\number\year,
          \thehours.\ifnum\theminutes<10{0}\fi\theminutes}

\author[S. Ponnusamy]{Saminathan Ponnusamy
}
\address{S. Ponnusamy, Stat-Math Unit,
Indian Statistical Institute (ISI), Chennai Centre,
110, Nelson Manickam Road,
Aminjikarai, Chennai, 600 029, India.
}
\email{samy@isichennai.res.in, samy@iitm.ac.in}

\author[K.-J. Wirths]{Karl-Joachim Wirths}
\address{K.-J. Wirths, Institut f\"ur Analysis und Algebra, TU Braunschweig,
38106 Braunschweig, Germany.}
\email{kjwirths@tu-bs.de}

\subjclass[2010]{30C45}
\keywords{Univalent function, subordination, Julia's lemma, Schwarz lemma
}
\thanks{The first author is on leave from IIT Madras.}

\begin{abstract}
For $0<\lambda \leq 1$, let ${\mathcal U}(\lambda)$ denote the family of functions $f(z)=z+\sum_{n=2}^{\infty}a_nz^n$
analytic in the unit disk $\ID$ satisfying the condition $\left |\left (\frac{z}{f(z)}\right )^{2}f'(z)-1\right |<\lambda $
in $\ID$.  Although functions in this family are known to be univalent in $\ID$, the coefficient conjecture about $a_n$
for $n\geq 5$ remains an open problem. In this article, we shall
first present a non-sharp bound for $|a_n|$.  Some members of the family ${\mathcal U}(\lambda)$ are given by
$$ \frac{z}{f(z)}=1-(1+\lambda)\phi(z) + \lambda (\phi(z))^2
$$
with $\phi(z)=e^{i\theta}z$,  that solve many extremal problems
in ${\mathcal U}(\lambda)$. Secondly, we shall consider the following question: Do there exist  functions
$\phi$ analytic in $\ID$ with $|\phi (z)|<1$ that are not of the form $\phi(z)=e^{i\theta}z$
for which the corresponding functions $f$ of the above form are members of the family ${\mathcal U}(\lambda)$?
Finally, we shall solve the second coefficient ($a_2$) problem in an explicit form for $f\in {\mathcal U}(\lambda)$
of the form
$$f(z) =\frac{z}{1-a_2z+\lambda z\int_0^z\omega(t)\,dt},
$$
where $\omega$ is analytic in $\ID$ such that $|\omega(z)|\leq 1$ and $\omega(0)=a$, where $a\in \overline{\ID}$.
\end{abstract}


\maketitle
\pagestyle{myheadings}
\markboth{S. Ponnusamy and K.-J. Wirths}{Coefficient problems  on the class $U(\lambda)$}
\cc

We denote the unit disk by $\ID=\{z\in \IC:\, |z| < 1\}$, and let
${\mathcal H}$ be the linear space of analytic functions defined on $\ID$ endowed with the topology of locally uniform convergence
and ${\mathcal A}=\{f\in {\mathcal H}:\, f(0)=f'(0)-1=0\}.$
The family  ${\mathcal S}$ of univalent functions from $\mathcal A$ and  many of its subfamilies, for which the image domains
have special geometric properties, have been investigated in detail. Among them are convex, starlike,
close-to-convex, spirallike and typically real mappings. For the general theory of univalent
functions we refer the reader to the books \cite{Duren:univ,Go,P}. However, the class ${\mathcal U}(\lambda)$ defined below
seems to have many interesting properties (cf. \cite{OPW-17b,PW-2017}).
For $0<\lambda \leq 1$, we consider the family
$${\mathcal U}(\lambda)=\{f\in {\mathcal A}:\,  \mbox{ $|U_f(z)| <\lambda$ in $\ID$}\},
$$
where
\be\label{OS-eq2}
U_f(z)=\left (\frac{z}{f(z)} \right )^{2}f'(z)-1=\frac{z}{f(z)} -z\left (\frac{z}{f(z)} \right )'-1, \quad z\in\ID.
\ee
Set  ${\mathcal U}:={\mathcal U}(1)$, and observe that ${\mathcal U}\subsetneq {\mathcal S}$ (see \cite{Aks58,AksAvh70}). Recently,  in \cite{PW-2017}, the
present authors have presented a simpler proof of it in a general setting.

More recently, a number of new and useful properties of the family  ${\mathcal U}(\lambda)$ were established in \cite{OPW,OPW-17a,OPW-17b}.
However, the coefficient problem for ${\mathcal U}(\lambda)$ remains open. This article
supplements the earlier investigations in this topic. See \cite{OPW,OPW-17a,OPW-17b}.

Let ${\mathcal B}=\{\omega\in {\mathcal H}:\, |\omega(z)|<1 ~\mbox{on $|z|<1$}\}$ and
${\mathcal B}_0=\{\omega\in {\mathcal B}:\, \omega (0)=0\}$. In addition, for $f, g \in {\mathcal H}$,
we use the symbol $f(z)\prec g(z)$, or in short $f\prec g$, to mean that there exists an $\omega\in {\mathcal B}_0$ such that  $f(z)=g(\omega (z))$.
We now recall the following results from \cite{OPW} which we need in the sequel.

\begin{Thm}\label{OPW11-thA}
Suppose that $f\in {\mathcal U}(\lambda)$ for some $\lambda \in (0,1]$ and $a_2=f''(0)/2$. Then we have the following:
\bee
\item[{\rm (a)}] If $|a_2|=1+\lambda$, then $f$ must be of the form
$$f(z)=
\frac{z}{(1+e^{i\phi}z)(1+\lambda e^{i\phi}z)}.
$$
\item[{\rm (b)}] $$ \frac{z}{f(z)}+ a_2z \prec 1+2\lambda z + \lambda z^2~ \mbox{ and }~
\frac{f(z)}{z}\prec \frac{1}{(1-z)(1-\lambda z)}, ~ z\in \ID .
$$
\eee
\end{Thm}

As an analog to Bieberbach conjecture for the univalent family ${\mathcal S}$ proved
by de Branges \cite{DeB1} (see also \cite{AvWir-09}),  the following conjecture was proposed in \cite{OPW}.

\bcon\label{conj1}
Suppose that $f\in {\mathcal U}(\lambda)$ for some $0<\lambda \leq 1$ and $f(z)=z+\sum_{n=2}^{\infty}a_nz^n$.  Then $|a_n|\leq \sum_{k=0}^{n-1}\lambda ^k$ for $n\geq 2$.
\econ

This conjecture has been verified for $n=2$ first in \cite{VY2013} and a simpler proof was given  in \cite{OPW}.
More recently,  in  \cite{OPW-17b},  Obradovi\'c et al. proved the conjecture for $n=3,4$ with an alternate proof for the case $n=2$,
but it remains open for all $n\ge 5$. Because ${\mathcal U}(1)\subsetneq {\mathcal S}$ and the Koebe function belongs to ${\mathcal U}(1)$,
this conjecture obviously holds for $\lambda =1$, in view of the de Branges theorem. Since no bound has been obtained for $|a_n|$
for $n\geq 5$, it seems useful to obtain a reasonable estimate. This attempt gives the
following theorem and at the same time the proof for the case  $\lambda =1$ does not require the use of
de Branges theorem that $|a_n|\leq n$ for $f\in {\mathcal S}$ with equality for the Koebe function and its rotation.

\bthm\label{OPW11-th2}
Let $f(z)=z+\sum_{n=2}^{\infty}a_nz^n$ belong to ${\mathcal U}(\lambda)$ for some $0<\lambda \leq 1$. Then
$$|a_{n}|\leq  1+\lambda \sqrt{n-1}\sqrt{\sum_{k=0}^{n-2}\lambda^{2k}},~ \mbox{ for $n\geq 2$}.
$$
\ethm
\bpf
Let $f\in {\mathcal U}(\lambda)$. Then the second subordination relation in Theorem \Ref{OPW11-thA}(b) shows that
$$\frac{f(z)}{z}\prec \frac{1}{1-\lambda z}\,\frac{1}{1- z}=f_1(z)\,f_2(z), \quad z\in\ID.
$$
Note that for
$$g_1(z)=\sum_{n=0}^{\infty}b_nz^n \prec f_1(z)=\frac{1}{1-\lambda z} ~\mbox{ and }~g_2(z)=\sum_{n=0}^{\infty}c_nz^n \prec f_2(z)=\frac{1}{1-z},
$$
where $b_0=c_0=1$, Rogosinski's theorems \cite{Rogo43} (see also \cite[Theorems 6.2 and 6.4]{Duren:univ}) give that
\begin{equation}\label{f1}
\sum_{k=1}^n|b_k|^2\leq  \sum_{k=1}^n\lambda^{2k} ~\mbox{ and }~ |c_n|\leq  1 ~\mbox{ for $n\geq 1$}.
\end{equation}
Moreover, the relation $\frac{f(z)}{z}=g_1(z)g_2(z)$ gives
$$a_{n+1}=\sum_{k=0}^nb_kc_{n-k}.
$$
Consequently, by \eqref{f1}, it follows from the classical Cauchy-Schwarz inequality that
$$|a_{n+1}|\leq 1+\sum_{k=1}^n|b_k|\leq 1+\sqrt{n}\sqrt{\sum_{k=1}^n|b_k|^2} \leq 1+\sqrt{n}\sqrt{\sum_{k=1}^n\lambda^{2k}} ,
$$
which implies the desired assertion.
\epf

Suppose that $f\in{\mathcal U}(\lambda)$. Then the second subordination relation in Theorem \Ref{OPW11-thA}(b) shows that
there exists a function $\phi \in \mathcal{B}_0$ such that
\be\label{OPW11-eq1}
\frac{z}{f(z)}=1-(1+\lambda)\phi(z) + \lambda (\phi(z))^2,\quad z\in\mathbb{D}.
\ee
From Theorem \Ref{OPW11-thA}(a), we see that  there is  a member in the family ${\mathcal U}(\lambda)$ in the above form with
$\phi(z)=e^{i\theta}z$. In this type of functions, we have $|a_2|=1+\lambda$. A natural question is whether there exist  functions
$\phi \in \mathcal{B}_0$ that are not of the form $\phi(z)=e^{i\theta}z$ of the above type for which the corresponding $f$ of the form
\eqref{OPW11-eq1} belongs to ${\mathcal U}(\lambda)$. In order to prove the next result, we need the classical Julia lemma
which is often quoted as Jack's lemma \cite[Lemma 1]{Jack} or Clunie-Jack's lemma \cite{clunie} although this fact
was known much before the work of Jack. See the article of Boas \cite{boas10} for historical commentary and the application of Julia's lemma.

\begin{Lem} {\rm (Julia's lemma)}\label{JMM-lemma}
Let $|z_0|<1$ and $r_0=|z_0|.$  Let $f(z)=\sum_{k=n}^{\infty}a_kz^k$ be continuous on
$|z|\leq r_0$ and analytic on $\{z:\,|z|<r_0\}\cup \{z_0\}$ with
$f(z)\not\equiv  0$ and $n\geq 1.$ If  $\ds |f(z_0)|=\max _{|z| \leq r_0}  |f(z)| ,$
then $z_0f'(z_0)/f(z_0)$ is real number and   $z_0f'(z_0)/f(z_0)\, \geq n.$
\end{Lem}


\bthm\label{OPW11-th1}
Suppose that $\phi \in \mathcal{B}_0$ that are not of the form $\phi(z)=e^{i\theta}z$ of the above type \eqref{OPW11-eq1}
such that there exists a $\theta_0$ with $\phi( e^{i\theta_0})=-1$. In addition we let $\phi$ be analytic on the closed unit disk
$\overline{\ID}$. Then $f$ expressed by the relation \eqref{OPW11-eq1} cannot be a member of the family ${\mathcal U}(\lambda)$.
\ethm
\bpf
We observe that  $f\in {\mathcal U}(\lambda)$ if and only if
$$\left |\frac{z}{f(z)} -z\left (\frac{z}{f(z)} \right )'-1\right | <\lambda,\quad z\in\mathbb{D},
$$
which according to \eqref{OS-eq2} and \eqref{OPW11-eq1} implies that there exists a function $\phi \in \mathcal{B}_0$ such that

\be\label{OPW11-eq2}
L(\phi )(z)
=\left |-(1+\lambda)(\phi(z)-z\phi'(z))+\lambda \phi(z)(\phi(z)-2z\phi'(z))\right |<\lambda, ~ z\in\mathbb{D}.
\ee
Note that we consider analytic functions $\phi$ in $\overline{\ID}$ that are not of the form $\phi(z)=e^{i\theta}z$ of the above type such that there
exists a $\theta_0$ with $\phi( e^{i\theta_0})=-1$. Examples of such functions are the Blaschke products with the above exception.
From Julia's lemma with $n=1$, we know that
$$\frac{z_0\phi'(z_0)}{\phi(z_0)} =m(\theta _0) \geq 1, \quad z_0=e^{i\theta _0}.
$$
If we let $\phi(z)=z\psi(z)$,  then we see that $\psi(\mathbb{D})\subset \overline{\mathbb{D}}$ and $\psi(e^{i\theta_0})) = -e^{-i\theta_0} $.
Now, we assume that $m(\theta_0)=1$. Since
$$\frac{z\phi'(z)}{\phi(z)}=1+\frac{z\psi'(z)}{\psi(z)},
$$
this means that $\psi'( e^{i\theta_0})=0$. But then an angle with width $\pi$ and vertex  $e^{i\theta_0}$ would be mapped
by $\psi$ onto an angle with width $2\pi$ or more and a vertex $-e^{-i\theta_0}  $. This contradicts the fact that
$\psi(\mathbb{D})\subset \overline{\mathbb{D}}$. Hence, $m(\theta_0)>1$. From the above we get
$$ e^{i\theta_0}\phi'(e^{i\theta_0}) = -m(\theta_0),
$$
and therefore,
\beqq
L(\phi )(z_0) &=&|-(1+\lambda)(\phi(z_0)-z_0\phi'(z_0))+\lambda \phi(z_0)(\phi(z_0)-2z_0\phi'(z_0))|\\
& = &\lambda + (1+3\lambda)(m(\theta_0)-1)
\eeqq
which shows that $L(\phi )(z_0) >\lambda.$ This contradicts \eqref{OPW11-eq2} and hence, $f$ cannot be a
member of the family ${\mathcal U}(\lambda)$. The proof is complete.
\epf

In \cite[Theorem 5]{OPW}, under a mild restriction on $f\in {\mathcal U}(\lambda)$, the region of variability of $a_2$ is
established as in the following form.

\begin{Thm}\label{OPW7-th2a}
Let $f\in {\mathcal U}(\lambda)$ for some $0<\lambda \leq 1$, and such that
\be\label{eq4c}
\frac{z}{f(z)}\neq (1-\lambda )(1+z),\,\,z\in \ID.
\ee
Then, we have
\be \label{eq4d}
\frac{z}{f(z)}-(1-\lambda )z\prec 1+2\lambda z+\lambda z^{2}
\ee
and the estimate $ |a_2-  (1-\lambda )|\leq 2\lambda $ holds. In particular, $ |a_2| \leq 1+\lambda $ and the estimate is sharp as the
function $f_{\lambda}(z) =z/((1+\lambda z)(1+z))$ shows.
\end{Thm}

Certainly, it was not unnatural to raise the question whether the condition \eqref{eq4c} is necessary for a function $f$
to belong to the family $\mathcal{U}(\lambda)$. This question was indeed raised in \cite{OPW}. In the next result, we show that the condition
\eqref{eq4c} cannot be removed from Theorem \Ref{OPW7-th2a}. Before, we present the proof, it is worth recalling from \cite{OPW}
that if $f\in {\mathcal U}(\lambda )$, then for each $R\in (0,1)$, the function
$f_R(z)=R^{-1}f(Rz)$ also belongs to ${\mathcal U}(\lambda )$.

 \bthm\label{OPW11-th4}
 Let $f(z)= z/((1-z)(1-\lambda z))$ and for a fixed $R\in (0,1)$, let $f_R(z)=R^{-1}f(Rz)$.
 Then we have
 \bee
\item[{\rm (a)}] For $0<\lambda \leq 1/2$ there exists, for any $R\in (0,1)$, an $r\in (0,1)$ such that $F(R,r)=0$, where
\be\label{OPW11-eq3}
 F(R,r)=\frac{r}{f_R(r)}-(1-\lambda)(1+r).
\ee
\item[{\rm (b)}] For $1/2<\lambda <1$ there exists, for any
$$1>R>\frac{1+\lambda\,-\,\sqrt{(1-\lambda)(1+7\lambda)}}{2\lambda},
$$
an $r\in (0,1)$ such that $F(R,r)=0$.
\eee
\ethm
\bpf
We consider $F(R,r)$ given by \eqref{OPW11-eq3} and observe that
$$F(R,r) =\lambda R^2r^2 - r[R(1+\lambda)+1-\lambda] + \lambda.
$$
We see that in the cases indicated in the statement of the theorem $F(R,0)=\lambda >0$ and $F(R,1) < 0$. Indeed
$$F(R,1)=\lambda R^2 - R(1+\lambda)+2\lambda-1=-R[(1-R)\lambda +1]  - (1-2\lambda)
$$
which is less than zero for any $R\in (0,1)$ and for $0<\lambda \leq 1/2$. Similarly, for the case  $1/2<\lambda <1$,
one can compute the roots  of the equation $F(R,1)=0$ and obtain the desired conclusion. This proves the assertion
of Theorem \ref{OPW11-th4}.
\epf

Because of the characterization of functions in $ {\mathcal U}(\lambda)$ via functions in $\mathcal{B}$, the following
result is of independent interest. As pointed out in the introduction, it is known that if $f\in {\mathcal U}(\lambda)$, then
$|a_2|\leq 1+\lambda$ with equality for $f(z)=z/[(1-z)(1-\lambda z)]$ and its rotation.

\bthm\label{OPW11-th3}
Let $f\in {\mathcal U}(\lambda)$, $\lambda\in (0,1),$  have the form
\be\label{OPW11-eq1a}
f(z)= z+\sum_{n=2}^{\infty}a_nz^n =\frac{z}{1-a_2z+\lambda z\int_0^z\omega(t)\,dt}
\ee
for some $\omega \in \mathcal{B}$ such that $\omega(0)=a\in \mathbb{D}$ and $v(x)$ be defined by
$$v(x)=\int_0^1 \frac{x+t}{1+xt}\,dt = \frac{1}{x}-\frac{1-x^2}{x^2}\log(1+x) <1 ~\mbox{ for $0<x<1$},
$$
and $v(0)=\lim_{x\ra 0^{+}}v(x)=1/2$. Then
$|a_2| \leq 1+\lambda v(|a|).$
The result is sharp.
\ethm
\bpf Recall the fact that $f(z)=z+\sum_{n=2}^{\infty}a_nz^n \in {\mathcal U}(\lambda)$ if and only if
\be\label{OPW7-eq3}
 \frac{z}{f(z)}=1-a_2z+\lambda z\int_0^z \omega(t)\,dt \,\neq 0, \quad z\in \ID,
\ee
where $\omega\in{\mathcal B}$. By assumption  $\omega (0)=a \in \mathbb{D}$. As in the proof of \cite[Theorem 1]{OPW}, assume on the contrary that
\be\label{OPW7-eq4}
|a_2|=\frac{1+\lambda v(|a|)}{r}, \quad r\in (0,1),
\ee
and consider the function $F$ defined by
$$F(z)= \frac{1}{a_2}\left [1+\lambda z\int_0^z\omega(t)\,dt\right ] , \quad z\in\ID.
$$
Then, according to the Schwarz-Pick lemma applied to $\omega\in{\mathcal B}$, we can easily obtain that
$$|\omega(z)|\leq \frac{|a|+|z|}{1+|az|},\quad z\in \mathbb{D},
$$
and thus, as in the proof of \cite[Theorem 2]{OPW}, it follows that
$$\left|\int_0^z \omega(t)\,dt\right|\leq v(|a|)<1, \quad z\in \mathbb{D},
$$
where $v(x)$ is defined as in the statement.  Consequently, for $|z|\leq r$, we get by \eqref{OPW7-eq4}
$$|F(z)|\leq  \frac{1}{|a_2|}\left [1+\lambda | z| \left |\int_0^z\omega(t)\,dt\right |\right ]  \leq \frac{1+r\lambda v(|a|)}{|a_2|}
= \frac{(1+r\lambda v(|a|))r}{1+\lambda v(|a|)} <r.
$$
Hence $F$ is a mapping of the closed disk  $\overline{\ID}_r$ into itself, where $\ID_r =\{z:\, |z|< r\}$.
Secondly, we have for $z_1$ and $z_2$ in $\overline{\ID}_r$,
\beqq
|F(z_1)-F(z_2)| &= &  \frac{\lambda r}{1+\lambda v(|a|)}\left|z_1\int_0^{z_1} \omega(t)\,dt  +(-z_1+z_1-z_2)\int_0^{z_2} \omega(t)\,dt\right|\\
&\leq & \frac{\lambda r}{1+\lambda v(|a|)}\left(|z_1|\left|\int_{z_2}^{z_1}\omega(t)\,dt\right| + |z_1-z_2|\left|\int_0^{z_2}\omega(t)\,dt\right|\right)\\
&\leq & \frac{\lambda r}{1+\lambda v(|a|)}\,(|z_1|+v(|a|))|z_1-z_2|\\
&\leq & \frac{\lambda r(r+v(|a|))}{1+\lambda v(|a|)}\,|z_1-z_2|\\
& <& r |z_1-z_2|.
\eeqq
Therefore, $F$ is a contraction of the disk $\overline{\ID}_r$ and according to Banach's fixed point theorem,
$F$ has a fixed point in $\overline{\ID}_r$. This implies  that there exists a
$z_0\in \ID_r$ such that $F(z_0)=z_0$ which contradicts \eqref{OPW7-eq3} at $z_0\in \ID$ (and thus, \eqref{OPW7-eq4} is not true for any $r\in (0,1)$).
Hence, we must have $|a_2| \leq 1+\lambda v(|a|)$ for $f\in {\mathcal U}(\lambda)$.

To prove that the second coefficient inequality is sharp, we consider
\be\label{OPW11-eq6}
\omega(z)=\frac{z+a}{1+az},\quad a\in(0,1),
\ee
and  we use that
$$v(a)=\int_0^1\omega(t)\,dt.
$$
Hence,
$$1-(1+\lambda v(a))z+\lambda z\int_0^z\omega(t)\,dt = 1-z -\lambda z\int_z^1\omega(t)\,dt\,=:G(z).
$$
We claim that $G(z)\neq 0$ in $\ID$. Since $G(0)=1$, we may assume on the contrary that there exists a $z\in \mathbb{D}\setminus\{0\}$ such that $G(z)=0$.
This is equivalent to
$$\frac{1}{\lambda z} =\frac{1}{1-z} \int_z^1\omega(t)\,dt.
$$
As
$$\left|\frac{1}{\lambda z}\right| > 1 ~\mbox{ and }~ \left|\frac{1}{1-z} \int_z^1\omega(t)\,dt\right| \leq 1,
$$
we have now proved that $G(z)\neq 0$ for $z\in \mathbb{D}$. In particular, this implies that the function $f$
defined by
$$f(z)=\frac{z}{1-(1+\lambda v(a))z+\lambda z\int_0^z\omega(t)\,dt}
$$
belongs to the family  ${\mathcal U}(\lambda)$, where $\omega$ is given by \eqref{OPW11-eq6}. This proves the sharpness.
\epf

Moreover, one can show that a similar sharp inequality is sharp for any $\omega$ as above.

Since  $|\int_{z_1}^{z_2}\omega(t)\,dt|\,\leq |z_1-z_2|$, the function $\int_0^z\omega(t)\,dt$ is uniformly continuous in the open unit disk.
Therefore this function can be extended continuously onto the closed unit disk. Hence, the function
$v(\omega) :=\max \{\left|\int_0^z\omega(t)\,dt\right|:\,  z\in \overline{\ID}\}$ is well defined.
Suppose that $f\in {\mathcal U}(\lambda)$ is given by
$$f(z)=\frac{z}{1-a_2 z+\lambda z\int_0^z\omega(t)\,dt}
$$
for some $0\leq \lambda <1$, where $\omega\in{\mathcal B}$. Then
\be\label{OPW11-eq4}
|a_2|\leq 1+\lambda v(\omega),
\ee
is  valid and this inequality is sharp.

In order to prove this inequality, we assume again that
$$|a_2|=\frac{1+\lambda v(\omega)}{r},\quad r\in (0,1),
$$
and do similar steps as in the proof of Theorem \ref{OPW11-th3}.
The inequality \eqref{OPW11-eq4} can be shown to be sharp in the following way: Consider
$$\tilde{\omega}(z)\,=\,e^{i\varphi}\omega\left(e^{i\theta}z\right),
$$
where $\varphi, \theta \in [0,2\pi)$ are chosen such that
$$v(\omega)\,=\,\int_0^1\tilde{\omega}(t)\,dt.
$$
Next, we may proceed as before to complete the proof. However, we omit the details to avoid a repetition of
the arguments.

A more detailed consideration of these cases can give more explicit bounds for $|a_2|$ as follows.

\bthm\label{OPW11-th5}
Let $f\in {\mathcal U}(\lambda)$, $\lambda\in (0,1),$  have the form \eqref{OPW11-eq1a}
for some analytic function $\omega$ such that $|\omega (z)|\leq 1$ and $\omega(0)=a\in \overline{\mathbb{D}}$. Let further
\beqq
B_a(z)= \left \{ \begin{array}{cl}\ds
\frac{1}{\overline{a}} -\frac{1-|a|^2}{\overline{a}^2z}\log (1+\overline{a}z) & \mbox{for $ a\in \mathbb{D}\setminus \{0\}$,}\\
 a& \mbox{for $|a|=1$},\\
 \ds \frac{z}{2} & \mbox{for $a=0$.} \end{array}\right .
\eeqq
Then
$$|a_2|\,\leq 1 + \lambda \max\{\left|B_a\left(e^{it}\right)\right|:\, t\in [0,2\pi]\}.
$$
The inequality is sharp.
\ethm
\bpf
The function $f$ considered here by \eqref{OPW11-eq1a} is a member of the class ${\mathcal U}(\lambda)$ if and only if $z/f(z)\neq 0$, which
is equivalent to
$$ a_2 \neq \frac{1}{z}\,+\,\lambda \int_0^z\omega(t)\,dt := C_{\omega}(z),\quad z\in \mathbb{D}.
$$
Using the above argument, it is clear that the function $C_{\omega}$ can be extended continuously onto the boundary $\partial\mathbb{D}$.
Moreover this function is univalent on $\overline{\mathbb{D}}$. The proof of this assertion is similar to the above arguments. Indeed if $C_{\omega}(z_1)=C_{\omega}(z_2)$ for some $z_1\neq z_2, z_1, z_2 \in \overline{\mathbb{D}}$, then
$$  \frac{\lambda}{z_1-z_2}\int_{z_1}^{z_2}\omega(t)\,dt \,=\,\frac{1}{z_1z_2}
$$
which is not possible. Thus, $C_{\omega}$ is univalent on $\overline{\mathbb{D}}$ and therefore, for each $\omega$,
the curve  $C_{\omega}\left(e^{i\theta}\right),\,\theta \in [0,2\pi]$, is a Jordan curve which divides the plane
into two components. Let us call the bounded closed component
$\overline{\mathbb{C}}\setminus C_{\omega}(\mathbb{D})\,=:\,A_2(\omega)$.
Obviously, the function $f$ is in the class  ${\mathcal U}(\lambda)$ if and only if
$$a_2 \,\in \,\bigcup_{\omega(0)=a}A_2(\omega).
$$
Now, we look at the curves $C_{\omega}\left(e^{i\theta}\right),\,\theta \in [0,2\pi]$.
Since $\omega(0)=a$, the modulus of the function
$$\frac{\omega(z)\,-\,a}{1\,-\,\overline{a}\omega(z)}
$$
is bounded by unity in the unit disk and this function vanishes at the origin. This means that $\omega$ can be represented in the form
$$\omega(z)\,=\,\frac{a\,+\,z\varphi(z)}{1\,+\,\overline{a}z\varphi(z)},
$$
where $\varphi$ is analytic in $\ID$ and $|\varphi (z)|\leq 1$ for $z\in\ID$. In other words, $\omega (z)$ is subordinated to
$(a+z)/(1+\overline{a}z)$, $z\in \ID$.
Since the function $(a+z)/(1+\overline{a}z)$ maps the unit disk onto the unit disk, a convex domain, we may use now a theorem proved by Hallenbeck and Ruscheweyh in \cite{HR} (compare with \cite[Theorem 3.1b]{MM}). In our case this theorem implies that the function
$$\frac{1}{z}\int_0^z\omega(t)\,dt$$
is subordinated to the function
$$\frac{1}{z}\int_0^z\frac{a\,+t}{1\,+\,\overline{a}t}\,dt\,=\,B_a(z).$$
Therefore, we get the representation
$$\int_0^z\omega(t)\,dt\,=\,\frac{1}{\varphi(z)}\int_0^{z\varphi(z)}\frac{a\,+t}{1\,+\,\overline{a}t}\,dt\,=\,zB_a(z\varphi(z)),
$$
where $\varphi$ is analytic in $\ID$ and $|\varphi (z)|\leq 1$ for $z\in\ID$. Since $B_a$ is analytic in the closed unit disk,
this representation together with the above considerations implies that
$$|a_2|\,\leq\, \sup_{z\in\mathbb{D},\theta \in [0,2\pi]}\left|e^{-i\theta}\,+\,\lambda e^{i\theta}B_a(z)\right|\,\leq\,1\,+\,\lambda \max\{\left|B_a\left(e^{it}\right)\right|: \,t\in [0,2\pi]\}.
$$

Now, we have to prove the sharpness of the inequality. To that end, let $t_0$ be chosen such that
$$\left|B_a\left(e^{it_0}\right)\right|\,=\,\max\{\left|B_a\left(e^{it}\right)\right|: \, t\in [0,2\pi]\},
~\mbox{ and }~ B_a\left(e^{it_0}\right)\,=\,e^{i\alpha}\left|B_a\left(e^{it_0}\right)\right|.
$$
We take $2\theta = -\alpha, \psi = t_0-\theta $,  consider the function
$$\omega(z) = \frac{a\,+\,ze^{i\psi}}{1\,+\overline{a}ze^{i\psi}},
$$
and let $ a_2\,=\,e^{-i\theta}\,+\,\lambda e^{i\theta}B_a\left(e^{it_0}\right).$
Then we have
$$|a_2|\,=\,\left|e^{-2i\theta}\,+\lambda \,e^{i\alpha}\left|B_a\left(e^{it_0}\right)\right|\, \right|\,=\,1\,+\,\lambda \left|B_a\left(e^{it_0}\right)\right|.$$
Further, we consider
$$D(z)\,=\,1\,-\,\left(e^{-i\theta}\,+\lambda e^{i\theta}B_a\left(e^{it_0}\right)\right)z\,+\,\lambda z \int_0^z\frac{a\,+\,te^{i\psi}}{1+\overline{a}te^{i\psi}}\,dt .
$$
It is easily seen that in our case
$$D(z)\,=\,1\,-\,\left(e^{-i\theta}\,+\lambda e^{i\theta}B_a\left(e^{it_0}\right)\right)z\,+\,\lambda  z^2
B_a\left(ze^{i\psi}\right ) ~\mbox{ and }~ D\left(e^{i\theta}\right)=0.
$$
The assumption, that there would exist a second zero $w$ of $D$ in the unit disk, leads to
$$\frac{1}{w}\,+\,\lambda \int_0^w\omega (t)\,dt\,=\,e^{-i\theta}\,+\,\lambda \int_0^{e^{i\theta}}\omega (t)\,dt,
$$
which is impossible, because the right hand side of the last relation is seen to be $a_2$.
This implies that the function $f(z)=z/D(z)$ is a member of the class  ${\mathcal U}(\lambda)$.
\epf

\end{document}